\documentclass{article}

\usepackage{amssymb}
\usepackage{amsfonts}
\usepackage{amsmath}
\usepackage{latexsym}
\usepackage{mathtools}
\usepackage{centernot}

\newtheorem{theorem}{Theorem}[section]
\newtheorem{corollary}[theorem]{Corollary}
\newtheorem{lemma}[theorem]{Lemma}
\newtheorem{proposition}[theorem]{Proposition}
\newtheorem{example}[theorem]{Example}

\newtheorem{remark}[theorem]{Remark}
\newtheorem{definition}[theorem]{Definition}

\newcommand{\FF}{\mathbb F}
\newcommand{\CC}{\mathbb C}

\newenvironment{matriz}[1]{\left[ \begin{array}{#1}}{\end{array} \right]}

\def\Gl{\mathop{\rm Gl}\nolimits}
\def\deg{\mathop{\rm deg }\nolimits}
\def\rank{\mathop{\rm rank}\nolimits}
\def\diag{\mathop{\rm diag }\nolimits}
\def\lcm{\mathop{\rm lcm }\nolimits}

\def\cpr{\mathop{ \; \rm \angle \; }\nolimits}
\def\sp{\mathop{\rm span}\nolimits}

\newcommand{\se}{\ensuremath{\stackrel{s.e.}{\sim}}}

\newcommand{\Ll}{\mathcal L}

\newcommand{\ba}{\mathbf a}
\newcommand{\bd}{\mathbf d}

\newcommand{\bc}{\mathbf c}

\newcommand{\bw}{\mathbf w}

\newcommand{\bb}{\mathbf b}
\newcommand{\br}{\mathbf r}
\newcommand{\bs}{\mathbf s}

\newcommand{\bp}{\mathbf p}
\newcommand{\bq}{\mathbf q}

\newcommand{\bba}{\overline{\mathbf a}}
\newcommand{\bbb}{\overline{\mathbf b}}

\newcommand{\ssc}{ \large \texttt{c}}
\newcommand{\ssd}{ \large \texttt{d}}

\newcommand{\ssu}{ \large \texttt{u}}
\newcommand{\ssv}{ \large \texttt{v}}

\title{
On the change of the Weyr characteristics of matrix pencils after rank-one perturbations\thanks{
Partially supported by ``Ministerio de Econom\'{\i}­a, Industria y Competitividad (MINECO)'' of Spain and ``Fondo Europeo de Desarrollo Regional (FEDER)'' of EU through grants MTM2017-83624-P.
}}

\author{Itziar Baraga\~na\footnote{Departamento de Ciencia de la Computaci\'on e I.A.,
 Universidad del Pa\'{\i}s Vasco, UPV/EHU, 
 Donostia-San Sebasti\'an, Spain, e-mail: itziar.baragana@ehu.eus. 
  },
Alicia Roca\footnote{
Departamento de  Matem\'atica Aplicada,
 IMM, Universitat Polit\`ecnica de Val\`encia,   Valencia, Spain, e-mail:
aroca@mat.upv.es.
}
}

\date{}
\begin{document}
\maketitle

\begin{abstract}

The change of the Kronecker structure of a  matrix pencil perturbed by another pencil of rank one has been characterized in terms of the homogeneous invariant factors and the chains of column and row minimal indices of the initial and the perturbed pencils. We obtain here a new characterization  in terms of the homogeneous invariant factors and  the  conjugate partitions of the corresponding chains of  column and row minimal indices of both pencils.
  
We also define  the generalized Weyr characteristic of an arbitrary matrix pencil and obtain bounds for the change of it when the pencil is perturbed by another pencil of rank one. The results improve known results on the problem, hold for arbitrary perturbation   pencils of rank one, and for any algebraically closed field.

\end{abstract}

\medskip


\medskip


\section{Introduction}
\label{intro}

Lot has been said about perturbations of matrix operators. In particular, changes of the Jordan structure of a square matrix or of the Weierstrass structure of a regular pencil under bounded rank perturbations have been studied, for example, in \cite{TeDo16, TeDoMo08, MoDo03, Sa02, Sa04} from a generic point of view,  and in \cite{BaRo18, BCU15,  Silva88_1, Th80, Za91} for general perturbations.
Results on perturbations of arbitrary pencils can be found in \cite{BaRo20_2, TeDo07, DoSt20}. See also the references therein.

Recently, there have been obtained bounds for the changes of the generalized Weyr characteristic of a complex square  matrix pencil (see Remark \ref{remdecreasing} below) perturbed by another pencil of the form
$w(su^*-v^*)$ (\cite[Theorem 7.8]{LeMaPhTrWiCAOT21}).  This has been done relating the Jordan chains of a square pencil with those of a linear relation. In this paper we  extend the notion of Jordan chain to possibly non square matrix pencils and  express the generalized Weyr characteristic of a pencil in terms of its Kronecker structure.

 Observe that complex pencils of rank one can also be of the form $(su-v)w^*$. 
 We also obtain bounds for the  change of the generalized Weyr characteristic of a matrix pencil perturbed by another arbitrary matrix pencil of rank one,  over an algebraically closed field, improving the bounds of (\cite[Theorem 7.8]{LeMaPhTrWiCAOT21}).
 
Important to the present work is the characterization  in \cite[Theorem 5.1]{BaRo20_2} of the changes of the Kronecker structure of a  pencil perturbed by another pencil of rank one (see also \cite{DoSt20}). 
We express here this characterization in terms of the conjugate  partitions of the corresponding chains of column and row minimal indices of the pencils involved (Theorem \ref{cmaintheogenpenc} below). Although  Theorem 5.1 in \cite{BaRo20_2} and the current result in Theorem \ref{cmaintheogenpenc} hold  for pencils over arbitrary fields, to simplify the analysis we only state here both (Theorem 5.1 in \cite{BaRo20_2} and  Theorem \ref{cmaintheogenpenc}) for algebraically closed  fields.

To achieve our results we have introduced two types of sequences of integers: partitions and chains. The partitions can be finite or infinite. We will identify two partitions if they only differ in the number of zeros. The chains are of fixed length. This distinction has to be kept in mind along the paper. 

The paper is organized as follows: in Section \ref{secpreliminaries} we present the notation and some preliminary results, including \cite[Theorem 5.1]{BaRo20_2}. In Section \ref{secjordanchains} we recall the definition of a Jordan chain, extend it to arbitrary matrix pencils, define the generalized Weyr characteristic  of a pencil, and express it in terms of the Kronecker invariants. Section \ref{sectraduccion}   is devoted to translate  \cite[Theorem 5.1]{BaRo20_2} into terms of the conjugate  partitions of the minimal indices of the pencils involved.  The main result of the paper is presented in Subsection \ref{subsecboundsmain}. It requires several technical results which appear in Subsection \ref{subsecboundstechnical}. 
Finally, the paper ends with a conclusion section.

\section{Preliminaries}
\label{secpreliminaries}

Let $\FF$ be an algebraically closed  field. $\FF[s]$ denotes the ring of polynomials in the indeterminate $s$ with coefficients in $\FF$ and  $\FF[s, t]$ denotes the ring of polynomials in two
variables $s, t$ with coefficients in $\FF$.
We denote by $\FF^{p\times q}$,  $\FF[s]^{p\times q}$, and $\FF[s, t]^{p\times q}$   the vector spaces  of $p\times q$ matrices with elements in $\FF$, $\FF[s]$ and $\FF[s, t]$, respectively.
$\Gl_p(\FF)$ will be the general linear group of invertible matrices
in $\FF^{p \times p}$.

\medskip

A {\em matrix pencil} is a  polynomial matrix $A(s)\in \FF[s]^{p\times q}$ of degree at most one ($A(s)=A_0+sA_1$, with $A_0, A_1\in \FF^{p\times q}$). The {\em normal rank} of $A(s)$, denoted by $\rank (A(s))$,  is the order of the largest non identically zero minor of $A(s)$, i.e., it is the rank of $A(s)$ considered as a matrix on the field of fractions of $\FF[s]$. 
The pencil is  {\em regular} if $p=q=\rank(A(s))$. Otherwise it is  {\em singular}.

Two matrix pencils $A(s)=A_0+sA_1, B(s)=B_0+sB_1\in \FF[s]^{p\times q}$ are {\em strictly equivalent} ($A(s)\se B(s)$) 
if there exist invertible matrices $P\in \Gl_p(\FF)$,   $Q\in \Gl_q(\FF)$ such that  $A(s)=PB(s)Q$. By $A(s)\centernot  \se B(s)$ we will understand that $A(s)$ and $B(s)$ are not strictly equivalent.

Given the pencil $A(s)=A_0+sA_1 \in \FF[s]^{p\times q}$ of $\rank A(s)=\rho$, a complete system of invariants for  the strict equivalence   is formed by a chain of homogeneous polynomials $\phi_1(s,t) \mid \dots \mid \phi_{\rho}(s,t),\ \phi_i(s,t) \in \FF[s,t], \ 1\leq i\leq \rho$, called the  {\em homogeneous invariant  factors}, and two collections of nonnegative integers $c_1\geq \dots \geq c_{q-\rho}$ and $u_1\geq \dots \geq u_{p-\rho}$, called  the {\em column and row minimal indices} of the pencil, respectively. In turn, the homogeneous invariant factors are determined by a chain of polynomials  $\alpha_1(s)\mid \ldots \mid \alpha_\rho(s)$ in $\FF[s]$, called the {\em invariant factors}, and a chain of polynomials  $t^{k_1}\mid \ldots \mid t^{k_\rho}$ in $\FF[t]$, called  the {\em infinite elementary divisors} (see  \cite[Ch. 2]{Friedland80} or \cite[Ch. 12]{Ga74}).
In fact, we can write 
$$\phi_i(s,t)=t^{k_i}t^{\deg(\alpha_i(s))}\alpha_i\left(\frac st\right), \ 1\leq i\leq \rho.$$ 
We will refer to the complete system of invariants for the strict equivalence as the {\em Kronecker structure} of the pencil.

The sum of the degrees of the homogeneous invariant factors  plus the sum of the minimal indices is equal to the rank of the pencil. 
Also, if $B(s)=A(s)^T$, then $A(s)$ and $B(s)$ share the homogeneous invariant factors and have interchanged minimal indices, i.e., the column (row) minimal indices of $B(s)$ are the row (column) minimal indices of $A(s)$.
If $A(s)\in \FF[s]^{p\times q}$ and $\rank (A(s))=p$ ($\rank (A(s))=q$), then $A(s)$ does  not have row (column) minimal indices. As a consequence, 
the invariants for the strict equivalence of  regular matrix pencils are reduced to  the homogeneous invariant factors.

\medskip

Denote by $\overline{\FF}= \FF\cup\{\infty\}$. The {\em spectrum} of
$A(s)=A_0+sA_1\in \FF[s]^{p \times q}$ is defined as
$
\Lambda(A(s))=\{\lambda\in \overline{\FF} \; : \; \rank (A(\lambda))< \rank (A(s))\},
$
where we agree that $A(\infty)=A_1$. The elements $\lambda\in \Lambda(A(s))$ are the {\em eigenvalues} of $A(s)$.

\medskip

Let  $\alpha_1(s)\mid \dots \mid\alpha_\rho(s)$ and $\phi_1(s, t)\mid \dots \mid\phi_\rho(s, t)$, $\rho=\rank A(s)$, be the invariant factors and the homogeneous invariant factors of $A(s)=A_0+sA_1\in \FF[s]^{p\times q}$, respectively. Factorizing both, the invariant factors and the homogeneous invariant factors, we can write

$$
\alpha_{\rho-i+1}(s)=\prod_{\lambda\in \Lambda(A(s))\setminus\{\infty\}}(s-\lambda)^{n_i(\lambda, A(s))}, \quad 1\leq i \leq \rho,
$$
$$
\phi_{\rho-i+1}(s,t)=t^{n_i(\infty, A(s))}\prod_{\lambda\in \Lambda(A(s))\setminus\{\infty\}}(s-\lambda t)^{n_i(\lambda, A(s))}, \quad 1\leq i \leq \rho.
$$
The integers
$
n_1(\lambda, A(s))\geq \dots \geq n_\rho(\lambda, A(s))
$
are called the {\em partial multiplicities} of $\lambda$ in $A(s)$. For $\lambda \in \overline{\FF}\setminus \Lambda(A(s))$ we take $n_1(\lambda, A(s))=\cdots=n_{\rho}(\lambda, A(s))=0$. We agree that $n_i(\lambda, A(s))=+\infty$ for $i<1$ and $n_i(\lambda, A(s))=0$ for $i>\rho$, for $\lambda \in \overline{\FF}$. We also agree that $\alpha_i(s)=\phi_i(s,t)=1$ for $i<1$ and  $\alpha_i(s)=\phi_i(s,t)=0$ for $i>\rho$.

A canonical form for  the strict equivalence of matrix pencils is the Kronecker canonical form.
It is a matrix pencil of the form
$${\scriptsize\begin{matriz}{ccccc}
J(s) & 0 & 0 & 0 & 0\\
0 & N(s) & 0 & 0 & 0\\
0 & 0 & L(s) & 0 & 0\\
0 & 0 & 0 & R(s) & 0 \\
0 & 0 & 0 & 0 & 0
\end{matriz}}\in \FF[s]^{p\times q},$$
where $J(s)$ is a diagonal of Jordan blocks (here $J_{\lambda_0, k}(s)$ corresponds to the elementary divisor $(s-\lambda_0)^k$),
\begin{equation} \label{jordanblock}
J_{\lambda_0, k}(s)=
{\scriptsize \begin{bmatrix}
 s-\lambda_0    & 1 &    &  & \\
 & \ddots &  \ddots &    \\
  &  &  \ddots &  1\\
 &  &  & s-\lambda_0\\
\end{bmatrix}}
\in \FF[s]^{k\times k}, 
\end{equation}
$N(s)=\diag(N_{k_1}(s), \ldots, N_{k_{\rho}}(s)),$
 where $t^{k_1} \mid \dots \mid t^{k_{\rho}}$ are the infinite elementary divisors and (the block will be empty if $k=0$)
\begin{equation}\label{dei}
N_k(s)={\scriptsize \begin{matriz}{cccccccc}
1 & s   &  & \\
  &  \ddots & \ddots &    \\
 &  &  \ddots & s \\
  &  &  & 1
\end{matriz}}\in \FF[s]^{k\times k},
\end{equation}
$L(s)=\diag(L_{c_1}(s), \ldots, L_{c_r}(s)),$ where $c_1\geq \dots\geq c_r>0=c_{r+1}=\dots = c_{q-\rho}$ are the column minimal indices and 
\begin{equation} \label{cmi}
L_k(s)= {\scriptsize \begin{matriz}{ccccc}
s & 1 &   &  \\
 &    \ddots & \ddots &    \\
 & &  s & 1
 \end{matriz}
}\in \FF[s]^{k \times (k+1)}, 
\end{equation}
$R(s)=\diag(R_{u_1}(s), \dots, R_{u_{r'}}(s)),$
 where $u_1\geq \dots\geq u_{r'}>0=u_{r'+1}=\dots=u_{p-\rho}$ are the row minimal indices and $R_k(s)=L_k(s)^T\in \FF[s]^{(k+1)\times k}$,
understanding that the non specified components are zero.
For details see  \cite[Ch. 2]{Friedland80} or \cite[Ch. 12]{Ga74} for infinite fields, and  \cite[Ch. 2]{Ro03} for arbitrary fields.

\medskip

The results in this paper are strongly linked to some properties of collections of nonnegative integers. We will distinguish two notions.

We call  {\em partition} of a positive integer $n$ to  a finite or infinite 
sequence of nonnegative integers $\ba = (a_1, a_2, \ldots)$, almost all being zero, such that $a_1\geq a_2 \geq  \ldots$ and $a_1+a_2+ \ldots=n$. 
The number of nonzero components of $\ba$
is the  {\em length} of $\ba$ (denoted  $\ell(\ba)$). Notice that $\ell(\ba)\leq n$. Given  a finite partition $\ba = (a_1, a_2, \ldots, a_n)$, if  necessary, we consider   $a_i=0$ if $i>n$.
We identify two partitions that differ only in the number of components equal to zero.
Given  two  partitions $\ba= (a_1, a_2, \ldots, a_n)$ and $\bb= (b_1, b_2, \ldots, b_n)$,  $\ba$ is {\em majorized} by $\bb$ (denoted $\ba \prec \bb$) if $\sum_{i=1}^k a_i \leq \sum_{i=1}^k b_i $
for $1 \leq k \leq n-1$ and $\sum_{i=1}^n a_i =\sum_{i=1}^n b_i $.

The {\em conjugate  partition} of $\ba$, $\bba = (\bar{a}_1, \bar{a}_2, \ldots)$, is defined as $\bar{a}_k := \#\{i : a_i \geq k\}, \ k\geq 1$.
We define $\ba \cup \bb$ to be the partition whose components are those of $\ba$ and $\bb$ arranged in decreasing  order (perhaps not strictly)  and $\ba +\bb$ will be the partition whose components are the sums of the corresponding components of $\ba$ and $\bb$.
The following properties are satisfied:
$\ba \prec \bb \Leftrightarrow \bbb \prec \bba$ and
$\overline{\ba \cup \bb}= \overline{\ba} + \overline{\bb}$.

We call  {\em chain}  to  a finite  sequence of  integers  $\ssc = (c_1, c_2, \ldots, c_m)$ such that $c_1\geq c_2 \geq  \ldots \geq c_m$. 
When necessary, we will consider   $c_i=+\infty$ if $i<1$ and $c_i=-\infty$ when $i>m$. We remark that a chain has a fixed number of integer components.
\begin{definition}[\mbox{$1$step-generalized majorization}]\hspace{-0.15cm}\footnote{Particular case of generalized majorization \cite[Definition 2]{DoStEJC10}}
Given two chains of  integers $\ssc = (c_1, \dots, c_m)$ and  $\ssd=(d_1, \dots, d_{m+1})$, 
we say that  $\ssd$ is {\em 1step-majorized} by $\ssc$   (denoted $\ssd \prec' \ssc$) if
$$c_{i}=d_{i+1}, \quad h \leq i \leq m, $$
where $h=\min\{i: c_i<d_i\}$ ($c_{m+1}=-\infty$).
\end{definition}

All along this paper, the  chains involved have nonnegative components. We define the  {\em conjugate of a chain} of nonnegative integers $\ssc=(c_1, \ldots, c_m)$ as the partition $\bar{\ssc} = (\bar{c}_1,  \ldots)$, where $\bar{c}_k := \#\{i : c_i \geq k\}, \ k\geq 1$.
When necessary, we will consider the term $\bar{c}_0= \#\{i : c_i \geq 0\}=m$. Notice that if $\bc$ is the partition  $\bc=(c_1, \dots, c_m, 0, \dots, )$ then $\bar{\ssc}=\bar{\bc} $ and $\bar{\bar{\ssc}}=\bc $. 

\medskip

In the next theorem the change of the Kronecker structure of a  pencil perturbed by another pencil of rank one is characterized. The result was independently obtained in  \cite{DoSt20} and \cite{BaRo20_2}.

\begin{theorem}[\mbox{\cite[Theorem 5.1, Corollary 5.4]{BaRo20_2}}]
\label{maintheogenpenc}
Let $A(s), B(s)\in \FF[s]^{p\times q}$ be matrix pencils such that  $A(s)\centernot \se B(s)$. 
Let $\rank A(s)=\rho_1$, $\rank B(s)=\rho_2$, let 
$\phi_1(s, t)\mid \dots \mid \phi_{\rho_1}(s, t)$,
$c_1 \geq \dots \geq c_{q-\rho_1}\geq  0$ and
$u_1 \geq \dots \geq u_{p-\rho_1}\geq  0$
 be  the homogeneous invariant factors, column minimal indices and row minimal indices, respectively, of $A(s)$, and let 
$\psi_1(s, t)\mid \dots \mid \psi_{\rho_2}(s, t)$, 
$d_1 \geq \dots \geq d_{q-\rho_2}\geq  0$ 
and 
$v_1 \geq \dots \geq v_{p-\rho_2}\geq  0$
 be the homogeneous invariant factors, column minimal indices and row minimal indices, respectively, of $B(s)$.
Let
$\ssc=(c_1, \dots, c_{q-\rho_1})$, $\ssd=(d_1, \dots, d_{q-\rho_2})$,
$\ssu=(u_1, \dots, u_{p-\rho_1})$, $\ssv=(v_1, \dots, v_{p-\rho_2})$ ($c_0=d_0=u_0=v_0=+\infty$) and
$\rho=\min\{\rho_1, \rho_2\}$.
\begin{enumerate}
\item \label{equal}
If $\ssc=\ssd$, $\ssu=\ssv$, then there exists a pencil $P(s)\in \FF[s]^{p \times q}$ of    $\rank (P(s))=1$ such that 
$A(s)+P(s)\se B(s)$ if and only if
\begin{equation}\label{eqintfihr1}
\psi_{i-1}(s, t)\mid\phi_i(s, t)\mid\psi_{i+1}(s, t), \quad 1\leq i \leq \rho.
\end{equation}
\item \label{rowequal}
If $\ssc\neq \ssd$, $\ssu=\ssv$,
 let
$$\ell=\max\{i \; :\; c_i\neq d_i\},$$
$$
f=\max\{i\in\{1, \dots, \ell\}\; : \; c_i<d_{i-1}\}, \ \
f'=\max\{i\in\{1, \dots, \ell\}\; : \; d_i<c_{i-1}\},
$$
$$
G=\rho-1-\sum_{i=1}^{\rho-1}\deg(\gcd(\phi_{i+1}(s, t),\psi_{i+1}(s, t)))- \sum_{i=1}^{p-\rho}u_i.
$$
Then, there exists a pencil $P(s)\in \FF[s]^{p \times q}$ of  $\rank (P(s))=1$ such that 
$A(s)+P(s)\se B(s)$ if and only if
 {\em (\ref{eqintfihr1})} holds and
\begin{equation*}\label{eqGg}
 G\leq \sum_{i=1}^{q-\rho}\min\{c_i, d_i\}+ \max\{c_f, d_{f'}\}.
\end{equation*}
\item  
\label{colequalcdrns}
If $\ssc= \ssd$, $\ssu\neq \ssv$,
let
 $$\bar \ell=\max\{i \; :\; u_i\neq v_i\},$$
$$
\bar f=\max\{i\in\{1, \dots, \bar \ell\}\; : \; u_i<v_{i-1}\}, \ 
\bar f'=\max\{i\in\{1, \dots, \bar \ell\}\; : \; v_i<u_{i-1}\},
$$
$$
\bar G=\rho-1-\sum_{i=1}^{\rho-1}\deg(\gcd(\phi_{i+1}(s, t),\psi_{i+1}(s, t)))-\sum_{i=1}^{q-\rho}c_i.
$$
Then, there exists a pencil $P(s)\in \FF[s]^{p \times q}$ of $\rank (P(s))=1$ such that 
$A(s)+P(s)\se B(s)$ if and only if
   {\em (\ref{eqintfihr1})} holds and
\begin{equation*}\label{eqbarGg}
 \bar G\leq \sum_{i=1}^{p-\rho}\min\{u_i, v_i\}+ \max\{u_{\bar f}, v_{\bar f'}\}.
\end{equation*}

\item \label{inequal}
If $\ssc\neq \ssd$, $\ssu\neq \ssv$,
then there exists a pencil $P(s)\in \FF[s]^{p \times q}$ of $\rank (P(s))=1$ such that 
$A(s)+P(s)\se B(s)$ if and only if  {\em (\ref{eqintfihr1})}  
and  one of the four following conditions hold:
\begin{enumerate}
\item[(a)]
\begin{equation}\label{cdsr}
{\ssc} \prec' {\ssd}, \quad {\ssu} \prec' {\ssv},
\end{equation}
\begin{equation}\label{sumpi1}
\sum_{i=1}^{\rho} \deg(\lcm(\phi_i(s,t), \psi_i(s,t)))\leq x\leq \sum_{i=1}^{\rho}\deg(\gcd(\phi_{i+1}(s,t), \psi_{i+1}(s,t))),
\end{equation}
where $x=\rho-\sum_{i=1}^{q-\rho_1}c_i-\sum_{i=1}^{p-\rho_2}v_i$.
\item[(b)]
\begin{equation}\label{dcrs}
{\ssd\prec'\ssc, \quad \ssv\prec'\ssu},
\end{equation}
\begin{equation}\label{sumpi1t}
\sum_{i=1}^{\rho} \deg(\lcm(\phi_i(s,t), \psi_i(s,t)))\leq y\leq \sum_{i=1}^{\rho}\deg(\gcd(\phi_{i+1}(s,t), \psi_{i+1}(s,t))),
\end{equation}
where $y=\rho-\sum_{i=1}^{q-\rho_2}d_i-\sum_{i=1}^{p-\rho_1}u_i$.
\item[(c)]
   {\em (\ref{cdsr})} and  {\em (\ref{sumpi1t})} hold.
 \item[(d)]
   {\em (\ref{dcrs})} and  {\em (\ref{sumpi1})} hold. 

\end{enumerate}

\end{enumerate}

\end{theorem}

\section{Jordan chains of matrix pencils}\label{secjordanchains}

The definition  of Jordan chain can be found in \cite[Section 1.4]{GoLaRo09} for square regular matrix polynomials  and  in  \cite[Definition 7.1]{LeMaPhTrWiCAOT21}
for square matrix pencils (both  regular and singular) over $\CC$. On the other hand, the notion of Weyr characteristic of an eigenvalue of a pencil was introduced in \cite{Ho90} as the conjugate partition of that of the partial multiplicities of the eigenvalue. 

The  target of this section is to generalize to arbitrary matrix pencils the notions of Jordan chain and   Weyr characteristic of  pencils, and  to  express the generalized  Weyr characteristic   in terms of the Kronecker invariants of the  pencil.
The results hold over arbitrary fields.
 
\begin{definition}\label{defjch}
Given a matrix pencil $A(s)=A_0+sA_1\in \FF[s]^{p\times q}$, an ordered set 
$(x_k, \dots, x_0)$ in $\FF^{q}$ is a {\em right Jordan chain} of $A(s)$ at $\lambda \in \overline{\FF}$, of length $k+1$,   if $x_0\neq 0$ and
$$
\begin{array}{rrrrrrr}
\lambda \in\FF:& A(\lambda)x_0=0, & A(\lambda)x_i=-A_1x_{i-1}, \quad 1 \leq i \leq k,\\
\lambda =\infty :& A_1x_0=0, & A_1x_i=-A_0x_{i-1}, \quad 1 \leq i \leq k.
\end{array}
$$
\end{definition}
The set
$(x_k, \dots, x_0)$ in $\FF^{p}$ is a {\em left Jordan chain} of $A(s)$ at $\lambda \in \overline{\FF}$, of length $k+1$,  if it is a right Jordan chain of $A(s)^T$ at $\lambda$.

In the sequel we will deal with right Jordan chains and we will refer to them just as Jordan chains, omitting the term  ``right''. It can be easily seen that the results obtained for right Jordan chains hold for left Jordan chains by transposition.

Following \cite{LeMaPhTrWiCAOT21}, we denote by
$\Ll^\ell_\lambda(A(s))$ the subspace spanned by the vectors of the Jordan chains at
  $\lambda \in \overline{\FF}$, up to length $\ell \geq 1$. We agree that $\Ll^0_\lambda(A(s))=\{0\}$.
If $\rank(A(\lambda))=q$, then there is no any Jordan chain at $\lambda$ for $A(s)$, and we take $\Ll^\ell_\lambda(A(s))=\{0\}$ for $\ell \geq 0$.
Observe that, for $\lambda \in \overline{\FF}$, 
$\Ll^{i-1}_\lambda(A(s))\subseteq \Ll^{i}_\lambda(A(s))$ for $i\geq 1$, and 
$\Ll^{i-1}_\lambda(A(s))=  \Ll^{i}_\lambda(A(s))$ for $i>q$.
Again as in \cite{LeMaPhTrWiCAOT21}, we denote by $w_i(\lambda, A(s))$ the dimension of the quotient space $\frac{\Ll^{i}_\lambda(A(s))}{\Ll^{i-1}_\lambda(A(s))}$, i.e. 
\begin{equation*} \label{omegas}
w_i(\lambda, A(s))=\dim \Ll^i_\lambda(A(s))-\dim \Ll^{i-1}_\lambda(A(s)), \quad 1\leq i\leq q.
\end{equation*}

The following theorem was obtained in  \cite{LeMaPhTrWiCAOT21}.

\begin{theorem}\cite[Theorem 7.8]{LeMaPhTrWiCAOT21}
\label{theoboundsLMPTW}
Given a matrix pencil $A(s)\in \CC[s]^{n\times n}$, let $P(s)\in \FF[s]^{n\times n}$ be a rank-one matrix pencil of the form 
$$
P(s)=w(su^*+v^*), \quad u,v,w \in \CC^n, \quad (u,v) \neq (0,0), \quad w\neq 0.
$$
For $\lambda \in \overline{\CC}=\CC \cup \{\infty\}$ and $i\geq 1$, the following statements hold
 \begin{enumerate}
\item[(i)]
If both pencils $A(s)$ and $A(s)+P(s)$ are regular, then
$$
\mid w_i(\lambda, A(s)+P(s))- w_i(\lambda, A(s))\mid \leq 1.
$$
\item[(ii)]
If  $A(s)$ is regular but $A(s)+P(s)$ is singular, then
$$
-i\leq  w_i(\lambda, A(s)+P(s))- w_i(\lambda, A(s))\leq 1.
$$
\item[(iii)]
If  $A(s)$ is singular and $A(s)+P(s)$ is regular, then
$$
-1\leq  w_i(\lambda, A(s)+P(s))- w_i(\lambda, A(s))\leq i.
$$
\item[(iv)]
If both $A(s)$ and $A(s)+P(s)$ are singular, then
$$
\mid  w_i(\lambda, A(s)+P(s))- w_i(\lambda, A(s))\mid\leq i.
$$
\end{enumerate}

\end{theorem}

In this paper we obtain bounds for  $w_i(\lambda, A(s)+P(s))- w_i(\lambda, A(s))$ for arbitrary matrix pencils $A(s)$ and arbitrary rank-one perturbations $P(s)$.

Since the relation between the Kronecker invariants of a pencil and those of the pencil obtained by a rank-one perturbation of it is known (see Theorem \ref{maintheogenpenc}), we are going to relate the values of  $w_i(\lambda, A(s))$ with the  Kronecker invariants of  $A(s)$. This relation can also be obtained from  \cite[Theorem 8.1]{LeMaPhTrWiCAOT21}, but to improve the readability of the paper,  a short proof of it is included here.

\medskip

First, in the next lemma we state that to compute $\dim\Ll^i_\lambda(A(s))$ we can substitute  $A(s)$ by a strictly equivalent pencil.
\begin{lemma}\label{propesteq}
Let $A(s), \bar{A}(s)\in \FF[s]^{p\times q}$
be matrix pencils and let 
$\lambda \in \overline{\FF}$.  If $A(s)$ and $\bar{A}(s)$ are strictly equivalent, 
then
$$\dim\Ll^i_\lambda(A(s))=
\dim\Ll^i_\lambda(\bar{A}(s)), \quad i\geq 0. 
$$
\end{lemma}
\noindent
    {\it Proof.}  The proof is straightforward.

    \hfill $\Box$

In the next proposition we  analyze  $\dim\Ll^i_\lambda(A(s))$  under certain structures of $A(s)$: when $A(s)$ has a diagonal decomposition, and when it has zero columns or zero rows. 

\begin{proposition}\label{propdiag} \
\begin{enumerate}
\item
Let $C(s)\in \FF[s]^{p_1\times q_1}$, $D(s)\in \FF[s]^{p_2\times q_2}$
be matrix pencils,
$A(s)=\begin{bsmallmatrix}C(s)&0\\0&D(s)\end{bsmallmatrix}$ and
$\lambda \in \overline{\FF}$.  
Then
$$
\dim\Ll^i_\lambda(A(s))=
\dim\Ll^i_\lambda(C(s))+\dim\Ll^i_\lambda(D(s)), \quad i\geq 0.
$$
\item
Let $C(s)\in \FF[s]^{p\times q_1}$
be a matrix pencil,
$A(s)=\begin{bsmallmatrix}C(s)&0\end{bsmallmatrix}\in \FF[s]^{p\times (q_1+q_2)}$ and
$\lambda \in \overline{\FF}$.  
Then
$$
\dim\Ll^i_\lambda(A(s))=
\dim\Ll^i_\lambda(C(s))+q_2, \quad i\geq 1.
$$
\item
Let $C(s)\in \FF[s]^{p_1\times q}$
be a matrix pencil,
$A(s)=\begin{bsmallmatrix}C(s)\\0\end{bsmallmatrix}\in \FF[s]^{(p_1+p_2)\times q}$ and
$\lambda \in \overline{\FF}$.  
Then
$$\dim\Ll^i_\lambda(A(s))=
\dim\Ll^i_\lambda(C(s)), \quad i\geq 0.
$$
\end{enumerate}
\end{proposition}

\noindent
{\it Proof.} 
\begin{enumerate}
\item
The set 
$(\begin{bsmallmatrix}x_k\\y_k\end{bsmallmatrix}, \dots, \begin{bsmallmatrix}x_0\\y_0\end{bsmallmatrix})$ is a Jordan chain   of $A(s)$  at $\lambda$ if
and only if $(x_k, \dots, x_0)$ is a Jordan chain   of $C(s)$  at $\lambda$ and
$(y_k, \dots, y_0)$ is a Jordan chain   of $D(s)$  at $\lambda$.
Therefore, for $i \geq 0$, 
$\Ll^i_\lambda(A(s))=(\Ll^i_\lambda(C(s))\times\{0\})\oplus (\{0\}\times\Ll^i_\lambda(D(s)).$
\item
The set
$(\begin{bsmallmatrix}x_k\\y_k\end{bsmallmatrix}, \dots, \begin{bsmallmatrix}x_0\\y_0\end{bsmallmatrix})$ is a Jordan chain   of $A(s)$  at $\lambda$ if
and only if $(x_k, \dots, x_0)$ is a Jordan chain   of $C(s)$  at $\lambda$. 
Therefore, for  $i \geq 1$, 
$\Ll^i_\lambda(A(s))=\Ll^i_\lambda(C(s))\times \FF^{q_2}.$
\item
The set $(x_k, \dots,x_0)$ is a Jordan chain  of $A(s)$  at $\lambda$ if
and only if $(x_k, \dots, x_0)$ is a Jordan chain  of $C(s)$  at $\lambda$. 
Therefore, for $i \geq 0$, 
$\Ll^i_\lambda(A(s))=\Ll^i_\lambda(C(s)).$
\end{enumerate}
\hfill $\Box$

The next step is to compute $w_i(\lambda, A(s))$ when $A(s)$ is a  block component of a pencil in Kronecker canonical form.

\begin{proposition}\label{propblocs}\
\begin{enumerate}
\item
Let $\lambda_0 \in \FF$ and $\lambda \in \overline{\FF}\setminus\{\lambda_0\}$. Let $J_{\lambda_0, k}(s)$ be the pencil defined in  {\em (\ref{jordanblock})}. Then 
$$\begin{array}{c}
w_i(\lambda_0, J_{\lambda_0, k}(s)) = 
1, \quad 1\leq i \leq k,  \vspace{0.2cm}\\
w_i(\lambda, J_{\lambda_0, k}(s)) =  0, \quad  1\leq i \leq k.
\end{array}
$$
\item Let  $\lambda \in \FF$. Let $N_k(s)$ be the pencil defined in  {\em (\ref{dei})}. Then
$$\begin{array}{c}
w_i(\infty, N_k(s)) = 1, \quad 1\leq i \leq k, \\
w_i(\lambda, N_k(s))=  0, \quad  1\leq i \leq k.
\end{array}
$$

\item   Let $\lambda \in \overline{\FF}$. Let $L_{k-1}(s)$ be the pencil defined in  {\em (\ref{cmi})}. Then
$$
w_i(\lambda, L_{k-1}(s)) =  
1, \quad  1\leq i \leq k.$$

\item
Let $\lambda \in \overline{\FF}$ and let $R_{k-1}(s)=L_{k-1}(s)^T$. Then
$$
w_i(\lambda, R_{k-1}(s))= 0, \quad  1\leq i \leq k-1.
$$

\end{enumerate}
\end{proposition}
{\it Proof.} 
Let $e_1, \dots, e_k$ be the columns of $I_k$.
\begin{enumerate}
\item
For $1\leq i \leq k$, 
$
\Ll_{\lambda_0}^i(J_{\lambda_0, k}(s))=\sp\{e_1,  \dots, e_i\}, 
$
hence
$
\dim\Ll^i_{\lambda_0}(J_{\lambda_0, k}(s))=i$, from where we obtain 
$
w_i(\lambda_0, J_{\lambda_0, k}(s))=1.
$
If 
$\lambda \in \FF\setminus\{\lambda_0\}$,  $\rank J_{\lambda_0,  k}(\lambda)=k$, hence
$
w_i(\lambda, J_{\lambda_0, k}(s))=0$ for $1\leq i\leq k$.
Analogously, for $1\leq i\leq k$, we obtain 
$w_i(\infty, J_{\lambda_0, k}(s))=0$.

\item
As in the previous case, for $1\leq i \leq k$,
$
\Ll_{\infty}^i(N_k(s))=\sp\{e_1, \dots, e_i\}$ and 
$\rank N_k(\lambda)=k$ for $\lambda \in \FF$.
\item
We have  
$
\Ll_\infty^i(L_{k-1}(s))=\sp\{e_k, e_{k-1}, \dots, e_{k-i+1}\}
$ for $1\leq i \leq k$. If 
$\lambda \in \FF$, then
$\Ll_\lambda^1(L_{k-1}(s))=\sp\left\{\begin{bsmallmatrix}1\\-\lambda \\\lambda^2\\\vdots \\(-\lambda)^{k-1}\end{bsmallmatrix}\right\}$, 
and for $2\leq i \leq k$,
$$
\Ll_\lambda^i(L_{k-1}(s))=\sp\left\{
\begin{bsmallmatrix}1\\-\lambda \\\lambda^2\\\vdots\\(-\lambda)^{i-1} \\\vdots \\(-\lambda)^{k-1}\end{bsmallmatrix}, 
\begin{bsmallmatrix}0\\-1\\x_{32}\\\vdots \\x_{i-1,2}\\x_{i2}\\\vdots \\x_{k2}\end{bsmallmatrix}, 
\begin{bsmallmatrix}0\\0\\1\\\vdots\\x_{i-1,3} \\x_{i3}\\\vdots \\x_{k3}\end{bsmallmatrix},  \dots, 
\begin{bsmallmatrix}0\\0\\0\\\vdots \\0\\(-1)^{i-1}\\\vdots \\x_{ki}\end{bsmallmatrix}\right\},
$$
where $x_{uj}$ are recursively defined as
$$
x_{uj}=-x_{u-1,j-1}-\lambda x_{u-1,j}, \quad 3\leq u \leq k, \; 2\leq j \leq u-1.
$$
Therefore, 
$\dim\Ll^i_{\lambda}(L_{k-1}(s))=i$ for $0\leq i \leq k$.

\item
For $\lambda \in \FF$, we have $\rank(R_{k-1}(\lambda))=k-1$. Hence, $w_i(\lambda, R_{k-1}(s))=0$ for $1\leq i\leq k$. Analogously, we obtain that
$w_i(\infty, R_{k-1}(s))=0$ for $1\leq i\leq k-1$.
\end{enumerate}
\hfill $\Box$

           \begin{corollary}
             \label{corblocs}
             With the notation of Proposition \ref{propblocs}:
\begin{enumerate}
\item
  Let $\lambda_0 \in \FF$ and $\lambda \in \overline{\FF}\setminus\{\lambda_0\}$.   Then 
  $$
  \begin{array}{l}
  \overline{(k)}=(w_1(\lambda_0, J_{\lambda_0, k}(s)), \dots, w_k(\lambda_0, J_{\lambda_0, k}(s))), \vspace{0.2cm}\\
  \overline{(0)}=(w_1(\lambda, J_{\lambda_0, k}(s)), \dots, w_k(\lambda, J_{\lambda_0, k}(s))).\end{array}
$$
\item Let  $\lambda \in \FF$. 
  Then
  $$
  \begin{array}{l}
    \overline{(k)}=(w_1(\infty, N_{ k}(s)), \dots, w_k(\infty, N_{ k}(s))), \vspace{0.2cm}\\ 
        \overline{(0)}=(w_1(\lambda, N_{k}(s)), \dots, w_k(\lambda, N_{k}(s))).\end{array}
$$
\item   Let $\lambda \in \overline{\FF}$.   Then
  $$
  \overline{(k)}=(w_1(\lambda, L_{ k-1}(s)), \dots, w_k(\lambda, L_{ k-1}(s))).$$
\item
  Let $\lambda \in \overline{\FF}$. 
  Then
  $$
  \overline{(0)}=(w_1(\lambda, R_{ k-1}(s)), \dots, w_{k-1}(\lambda, R_{ k-1}(s))).$$
\end{enumerate}
\end{corollary}

Next theorem relates the values of  $w_i(\lambda, A(s))$ with the Kronecker invariants of $A(s)$.

 \begin{theorem}\label{carwcanonical}
Let $A(s)\in \FF[s]^{p\times q}$ be a matrix pencil such that $\rank (A(s))= \rho$, and $\lambda \in \overline{\FF}$.   Let
$
n_1(\lambda, A(s))\geq \dots \geq n_\rho(\lambda, A(s))\geq 0
$
be the partial multiplicities of $\lambda$ in $A(s)$ and let
$c_1\geq \dots \geq c_{q-\rho}$ and $u_1\geq \dots \geq u_{p-\rho}$ be the column and row minimal indices of $A(s)$, respectively.
Let
$
(w^R_1(\lambda, A(s)),  \dots,  w^R_q(\lambda, A(s)) )
=\overline{(n_1(\lambda, A(s)),  \dots, n_\rho(\lambda, A(s))}
$,
$
(r_1, \dots,  r_q)=\overline{(c_1, \dots, c_{q-\rho})}
$,
and $r_0=q-\rho$.
Then 
$$
w_i(\lambda, A(s))=w^R_i(\lambda, A(s))+r_{i-1},\quad 1\leq i \leq q.
$$
\end{theorem}

{\it Proof.}
By Lemma \ref{propesteq}, we can assume that $A(s)$ is in Kronecker canonical form.
By Propositions \ref{propdiag} and \ref{propblocs}, for $\lambda \in \overline{\FF}$,
$$
\sum_{i=1}^{\rho}\overline{(n_i(\lambda, A(s))}+ \sum_{i=1}^{r_1}\overline{(c_i+1)}+(q-\rho-r_1, 0, \dots, 0)=
(w_1(\lambda, A(s)), \dots, w_q(\lambda, A(s))).
$$
Since
$$
\begin{array}{rl}
\sum_{i=1}^{\rho}\overline{(n_i(\lambda, A(s))}=&
\overline{ \cup_{i=1}^{\rho}(n_i(\lambda, A(s))}=\overline{(n_1(\lambda, A(s)),  \dots, n_\rho(\lambda, A(s))}\vspace{0.2cm}\\
\vspace{0.2cm}
= &(w^R_1(\lambda, A(s)),  \dots, w^R_q(\lambda, A(s))),
 \end{array}
 $$
and 
$
(q-\rho-r_1, 0, \dots, 0)=\sum_{i=r_1+1}^{q-\rho}1=\sum_{i=r_1+1}^{q-\rho}\overline{(c_i+1)},$
we obtain 
$$
\begin{array}{rl}
\sum_{i=1}^{r_1}\overline{(c_i+1)}+(q-\rho-r_1, 0, \dots, 0)=&\sum_{i=1}^{q-\rho}\overline{(c_i+1)}=
\overline{ \cup_{i=1}^{q-\rho}(c_i+1)}\vspace{0.2cm}\\
=
\overline{(c_1, \dots, c_{q-\rho})+(1, \stackrel{(q-\rho)}{\dots}, 1)}=& \overline{(c_1, \dots, c_{q-\rho})}\cup \overline{(1, \stackrel{(q-\rho)}{\dots}, 1)}
\\
=(r_1, \dots,  r_q)\cup
(q-\rho)=&(q-\rho, r_1, \dots, r_q).
\end{array}
$$
\hfill $\Box$

\begin{remark}\label{remdecreasing}\

\begin{itemize}
\item   Let $A(s)\in \FF[s]^{p\times q}$ be a matrix pencil and $\lambda \in \overline{\FF}$.   Then, by Theorem \ref{carwcanonical}, we have
  $$w_1(\lambda, A(s))\geq w_2(\lambda, A(s))\dots \geq  w_q(\lambda, A(s)). $$
We call the partition 
 \begin{equation*}
 \bw(\lambda, A(s))=(w_1(\lambda, A(s)), w_2(\lambda, A(s)), \cdots, w_q(\lambda, A(s))),
 \end{equation*}
the  {\em generalized Weyr characteristic} of the pencil $A(s)$ at $\lambda$.

\item Notice that, by Theorem \ref{carwcanonical},
$$\bw(\lambda, A(s))=(w^R_1(\lambda, A(s)),  \dots, w^R_q(\lambda, A(s)))+(r_0, r_1, \dots, r_{q-1}),\quad \lambda \in \overline{\FF},$$
where $(w^R_1(\lambda, A(s)),  \dots, w^R_q(\lambda, A(s)))$ is the Weyr characteristic at $\lambda$ of the regular part of $A(s)$ and $(r_0, r_1, \dots, r_{q-1} )$ is that of the column minimal indices block. Observe that the row minimal indices do not play any role in the calculation of $\bw(\lambda,A(s))$. 

In \cite{Ho90},  $(w^R_1(\lambda, A(s)),  \dots, w^R_q(\lambda, A(s)))$ is called the {\em Weyr characteristic} of the pencil for $\lambda \in \overline{\CC}$, and $(r_1, \dots, r_{q-1} )$ is called the {\em partition of the r-numbers} of the pencil.
Notice also that, if the pencil does not have column minimal indices, the Weyr characteristic and the generalized Weyr characteristic coincide. 

\medskip

In the next example we compute the Weyr and the generalized Weyr characteristics of a pencil.

 \begin{example}
   The spectrum of the pencil
 $$A(s)=
\begin{matriz}{ccc|ccc|c}
 1 & s & & &  & &  \\
  & 1& & & & & \\
 & & 1 & &  & &  \\
\hline
 & &  & s & 1 & &  \\
 & & &  & s & 1 &  \\
\hline
 & & & & & & 0 \\
\end{matriz}\in \CC[s]^{6\times 7}
$$ is $\Lambda(A(s))=\{\infty\}.$
The partial multiplicities of $\infty$ in 
$A(s)$ are $n_1(\infty, A(s))=2\geq n_2(\infty, A(s))=1$. The column and row minimal indices of $A(s)$ are $c_1=2\geq c_2=0$ and $u_1=0$, respectively.

\medskip

In \cite{Ho90}, the Weyr characteristic of $A(s)$ for  $\infty $ is  
$$(w_1^R (\infty, A(s)), \dots )=\overline{(2,1)}=(2,1, 0, \dots ).$$
 The partitions of the $r$-numbers and of the $s$-numbers of $A(s)$ are, respectively,
$$(r_1, \dots )=\overline{(2,0)}=(1,1, 0, \dots ), \quad (s_1, \dots )=\overline{(0)}=(0, \dots ).$$

The generalized Weyr characteristic of $A(s)$ at $\infty$ is
$$\bw (\infty, A(s))=(w_1^R (\infty, A(s)), \dots )+(r_0, r_1, \dots)$$
$$=
(2,1, 0, \dots )+(2,1,1, 0, \dots )=(4,2,1, 0,\dots).$$
\end{example}
\end{itemize}
\end{remark}

\section{Theorem \ref{maintheogenpenc} in terms of the conjugate  partitions}\label{sectraduccion}

The target of this section is to rewrite the characterizations stated in Theorem \ref{maintheogenpenc} in terms of the  conjugate partitions of the corresponding chains of column and row minimal indices of the pencils involved.
To achieve it we previously prove  some technical results.

We start with the introduction of  a new majorization between partitions of nonnegative integers.
 
     \begin{definition}\label{defmgc}
Given two partitions  of nonnegative integers $\br = (r_0, r_1, \dots)$ and  
$\bs=(s_0, s_1, \dots)$ 
such that
$r_0 \geq r_1\geq \dots $ and  $s_0 \geq s_1 \geq \dots $,
we say that  $\bs$ is {\em  conjugate majorized}  by $\br$ (denoted $\bs \cpr \br$) 
if $r_0=s_0+1$ and 
$$r_{i}=s_i+1, \quad 0 \leq i \leq g, $$
where $g=\max\{i: r_i>s_i\}$.
\end{definition}

     \begin{remark}
       Let  $\ba=(a_1, \dots )$ be a partition of nonnegative integers and let $(r_1, \dots)=\overline{(a_1, \dots )}$ be  its  conjugate partition.  In what follows we will frequently use the following properties:

\begin{itemize}
\item 
  $r_j=i,$ for $a_{i+1}<j \leq a_i,$ $i\geq 1.$
If $a_{i+1}=a_i$ there are no  $j$ such that  $r_j=i$.

\item For  $i\in \{1, \dots, \ell(\ba)\}$ we have $r_{a_i}\geq i$, and if  $j>a_i$ then $r_j<i$. Recall that  $\ell(\ba)$ is  the length of the partition $\ba$, i.e. the number of nonzero elements of $\ba$. In other words, $\ell(\ba)=r_1$.
\end{itemize}
\end{remark}

The proof of the next lemma is analogous to that of \cite[Lemma 3.2]{LiSt09}.

\begin{lemma}
\label{lemmaconj}
Let $(a_1, \dots)$ and $(b_1, \dots)$ be partitions of nonnegative integers. Let $\bp=(p_1, \dots)=\overline{(a_1, \dots)}$ and $\bq=(q_1, \dots)=\overline{(b_1, \dots)}$ be the  conjugate partitions.
Let $k\geq 0$ be an integer. 
Then,
\begin{equation*}\label{abk}
a_{j}\geq b_{j+k}, \quad \text{for every  }\  j\geq 1,
\end{equation*}
if and only if 
\begin{equation*}\label{ghk}
p_{j}\geq q_{j}-k, \quad  \text{for every  }\  j\geq 1.
\end{equation*}
\end{lemma}
\begin{lemma}\label{lemmagh}
Given two chains of  nonnegative  integers $\ssc = (c_1, \dots, c_{m+1})$ and  
$\ssd=(d_1, \dots, d_{m})$, let
$\br=(r_1, \dots)=\overline{(c_1, \dots, c_{m+1})}$, 
$\bs=(s_1, \dots)=\overline{(d_1, \dots, d_{m})}$
and $r_0=m+1=s_0+1$.
Let $$g=\max\{i: r_i>s_i\}, \quad h=\min\{i: d_i<c_i\}.$$ Then
\begin{equation}\label{eqgch}
  g=c_h,
  \end{equation}
and
\begin{equation}\label{eqsumrscd}
\sum_{j=1}^g(r_j-s_j-1)=\sum_{j=h}^m(c_{j+1}-d_j).
\end{equation}
\end{lemma}

{\it Proof.} 
As $d_h<c_h$, we have $s_{c_h}<h\leq r_{c_h}$, hence $g\geq c_h$. Observe that $g\leq \ell(\br)$.
If $\ell(\br)\geq i>c_h$, 
then $c_{r_i}\geq i >c_h$, hence $r_i<h$. By the definition of $h$, we have
$d_{r_i}\geq c_{r_i}\geq i$, from where $\#\{j\; : \; d_j\geq i\}\geq r_i$, i.e., $ s_i\geq r_i$. Therefore, (\ref{eqgch}) holds.
Then, 
$$
\sum_{j=1}^gr_j=\sum_{j=1}^{c_h}r_j=(m+1)c_{m+1}+\sum_{j=h}^{m}j(c_j-c_{j+1})=hc_h+\sum_{j=h}^{m}c_{j+1},
$$
and, bearing in mind that $d_h<c_h\leq c_{h-1}\leq d_{h-1}$,
$$
\sum_{j=1}^gs_j=\sum_{j=1}^{d_h}s_j+
\sum_{j=d_h+1}^{c_h}s_j=md_m+(c_h-d_h)(h-1)+\sum_{j=h}^{m-1}j(d_j-d_{j+1})$$$$=
hd_h+(c_h-d_h)(h-1)+\sum_{j=h}^{m-1}d_{j+1}
=(h-1)c_h+\sum_{j=h}^{m}d_{j}.
$$
Therefore,
$$
\sum_{j=1}^g(r_j-s_j-1)=
c_h-g+\sum_{j=h}^m(c_{j+1}-d_j)=\sum_{j=h}^m(c_{j+1}-d_j).
$$
\hfill $\Box$

  \begin{proposition}\label{propconj}
Given two chains of  nonnegative  integers $\ssc = (c_1, \dots, c_{m+1})$ and  
$\ssd=(d_1, \dots, d_{m})$, let
$(r_1, \dots)=\overline{(c_1, \dots, c_{m+1})}$,
$(s_1, \dots)=\overline{(d_1, \dots, d_{m})}$ be the  conjugate partitions, 
$r_0=m+1=s_0+1$, and
$\br=(r_0, r_1, \dots)$, $\bs=(s_0, s_1, \dots)$, 
Then $\ssc\prec'\ssd$ if and only if  $\bs \cpr \br$.
\end{proposition}

\noindent
{\it Proof}.
Let $g=\max\{i: r_i>s_i\}$ and $h=\min\{i: d_i<c_i\}.$ Then,
by Lemma \ref{lemmagh},  (\ref{eqgch}) and 
(\ref{eqsumrscd}) hold.
Moreover, $r_i\leq s_i<s_i+1$ for $i>g$ and $d_i\geq c_i\geq c_{i+1}$ for $1\leq i <h$.

Assume that  $\bc\prec'\bd$, then
$d_i\geq c_{i+1}$ for $1\leq i \leq m$ and $\sum_{j=h}^m(c_{j+1}-d_j)=0$.
From Lemma \ref{lemmaconj} we have
$r_i\leq s_{i}+1$ for $i\geq 1$, and then, from (\ref{eqsumrscd}) we derive that
$r_i= s_{i}+1$ for $1\leq i \leq g$, i.e., $\bs \cpr \br$.

Conversely, assume that  $\bs \cpr \br$,  then
$r_i\leq s_{i}+1$ for $i\geq 1$ and $\sum_{j=1}^g(r_j-s_j-1)=0$.
From Lemma \ref{lemmaconj} we have 
$d_i\geq c_{i+1}$ for $1\leq  i \leq m$, and then, from (\ref{eqsumrscd}) we derive that
$d_i= c_{i+1}$ for $h\leq i \leq m$, i.e., $\bc\prec'\bd$.

\hfill $\Box$

\begin{lemma}\label{lemmagh2}
Given two chains of  nonnegative  integers $\bc = (c_1, \dots, c_{m})$ and  
$\bd=(d_1, \dots, d_{m})$ $(d_0=c_0=+\infty)$ such that $\bc\neq \bd$, let
$$\ell=\max\{i \; :\; c_i\neq d_i\},$$
$$f=\max\{i\in\{1, \dots, \ell\}\; : \; c_i<d_{i-1}\},  \ \ f'=\max\{i\in\{1, \dots, \ell\}\; : \; d_i<c_{i-1}\}.
$$
Let
$(r_1, \dots)=\overline{(c_1, \dots, c_{m})}$, 
$(s_1, \dots)=\overline{(d_1, \dots, d_{m})}$, 
$r_0=s_0=m$,  
$$x=\min\{i \; :\; r_i\neq s_i\},$$
$$
e=\min\{i\geq x-1 : \; s_{i+1}\geq r_{i+1}\}, \quad e'=\min\{i\geq x-1: \; r_{i+1}\geq s_{i+1}\}.
$$
Then
$$
  e=c_f, \quad e'=d_{f'}.
  $$
\end{lemma}
{\it Proof.} 
Assume that $d_{\ell}<c_{\ell}$ (if $c_\ell< d_\ell$ the proof is analogous). As $d_{\ell}<c_{\ell}\leq c_{\ell-1}$, we have  $f'=\ell$.  

From the definition of $\ell$, $c_i=d_i, \ \ell+1\leq i\leq  m$, then 
$r_i=s_i, \ 1\leq i\leq d_{\ell}$. 
From $d_{\ell}<c_{\ell}$ we also derive that $s_{d_{\ell}+1}<r_{d_{\ell}+1}$, hence $x=d_{\ell}+1$ and, from the  definition of $e'$,  $e'=d_{\ell}$.

From the definition of $f$, $c_f<d_{f-1}$. If $f<\ell$, then $d_f\leq c_{f+1}\leq c_f<d_{f-1}$, and if $f=\ell$, then $d_f=d_{\ell}< c_{\ell}=c_f<d_{f-1}$. Hence, $r_{c_f+1}\leq f-1=s_{c_f+1}$ and $x-1=d_{\ell}\leq e\leq c_f$. Moreover, by definition, $c_i\geq d_{i-1}$ for $f+1\leq i\leq \ell$.

We prove next that  for $d_{\ell}\leq i\leq c_f-1$, $r_{i+1}>s_{i+1}$ holds. 
Let $i$ be such that $c_f\geq i> d_{\ell}$. If there exists $j\in \{f, \ldots, \ell\}$ such that $c_j\geq i>c_{j+1}\geq d_j$, then $r_i=j>j-1\geq s_i$. Otherwise, $c_{\ell}\geq i>d_{\ell}$, then $r_i> s_i$. Hence $e=c_f$. 
\hfill $\Box$
\begin{lemma}\label{lemmapartmin}
Given two chains of  nonnegative  integers $\bc = (c_1, \dots, c_{m})$ and  
$\bd=(d_1, \dots, d_{m})$, let
$x_i=\min\{c_i, d_i\}$, $1\leq i \leq m.$
Let
$(r_1, \dots)=\overline{(c_1, \dots, c_{m})}$, 
$(s_1, \dots)=\overline{(d_1, \dots, d_{m})}$, 
and
$y_i=\min\{r_i, s_i\}$, $i\geq 1.$
Then
$$
(y_1, \dots)=\overline{(x_1, \dots, x_{m})}.
  $$
\end{lemma}
    {\it Proof.}
 For $i\geq 1$, let 
    $
    C_i=\{j\; : \: c_j\geq i\}, \
    D_i=\{j\; : \: d_j\geq i\}, \
    X_i=\{j\; : \: x_j\geq i\}.
    $
    Then
    $
    X_i=C_i\cap D_i, \ r_i=\#C_i, \ s_i=\#D_i, \ i\geq 1.
    $
    We must prove that $y_i=\#X_i$ for $i\geq 1$.
    Let $i\geq 1$. If $y_i=r_i$, then $C_i\subseteq D_i$ and $X_i=C_i$, hence $\#X_i=r_i=y_i$.
Analogously, if $y_i=s_i$, then 
$\#X_i=s_i=y_i$. 

\hfill $\Box$

The result of the next  theorem is the target of the section.

 \begin{theorem}
\label{cmaintheogenpenc}
Let $A(s), B(s)\in \FF[s]^{p\times q}$ be matrix pencils such that  $A(s)\centernot \se B(s)$. 
Let $\rank A(s)=\rho_1$, $\rank B(s)=\rho_2$, let 
$\phi_1(s, t)\mid \dots \mid \phi_{\rho_1}(s, t)$,
$c_1 \geq \dots \geq c_{q-\rho_1}\geq  0$ and
$u_1 \geq \dots \geq u_{p-\rho_1}\geq  0$
 be  the homogeneous invariant factors, column minimal indices and row minimal indices of $A(s)$, respectively,  and let 
$\psi_1(s, t)\mid \dots \mid \psi_{\rho_2}(s, t)$, 
$d_1 \geq \dots \geq d_{q-\rho_2}\geq  0$ 
and 
$v_1 \geq \dots \geq v_{p-\rho_2}\geq  0$
 be  the homogeneous invariant factors, column minimal indices and row minimal indices of $B(s)$, respectively.

Let $\rho=\min\{\rho_1, \rho_2\}$, $\rho'=\max\{\rho_1, \rho_2\}$,
$\ssc=(c_1, \dots, c_{q-\rho_1})$, $\ssd=(d_1, \dots, d_{q-\rho_2})$,
$\ssu=(u_1, \dots, u_{p-\rho_1})$, $\ssv=(v_1, \dots, v_{p-\rho_2})$, 
$(r_1, \dots)=\overline{(c_1, \dots, c_{q-\rho_1})}$, 
$(s_1, \dots)=\overline{(d_1, \dots, d_{q-\rho_2})}$,
$(r'_1, \dots)=\overline{(u_1, \dots, u_{p-\rho_1})}$, 
$(s'_1, \dots)=\overline{(v_1, \dots, v_{p-\rho_2})}$,
$r_0=q-\rho_1$, $s_0=q-\rho_2$, $r'_0=p-\rho_1$, $s'_0=p-\rho_2$,
$\br=(r_0, r_1, \dots)$, 
$\bs=(s_0, s_1, \dots)$,
$\br'=(r'_0, r'_1, \dots)$, and
$\bs'=(s'_0, s'_1, \dots)$.

\begin{enumerate}
\item \label{cequal}
If $\br=\bs$, $\br'=\bs'$, then there exists a pencil $P(s)\in \FF[s]^{p \times q}$ of    $\rank (P(s))=1$ such that 
$A(s)+P(s)\se B(s)$ if and only if
 {\em (\ref{eqintfihr1})} holds.

\item \label{crowequal}
If $\br\neq \bs$, $\br'=\bs'$,
 let
\begin{equation} \label{T4.8-x}
x=\min\{i \; :\; r_i\neq s_i\},
\end{equation}
\begin{equation} \label{T4.8-ee'}
e=\min\{i\geq x-1 \; : \; s_{i+1}\geq r_{i+1}\},\quad 
e'=\min\{i\geq x-1\; : \; r_{i+1}\geq s_{i+1}\},
\end{equation} 
$$
G=\rho-1-\sum_{i=1}^{\rho-1}\deg(\gcd(\phi_{i+1}(s, t),\psi_{i+1}(s, t)))- \sum_{i=1}^{\rho}r'_i.
$$
Then, there exists a pencil $P(s)\in \FF[s]^{p \times q}$ of  $\rank (P(s))=1$ such that 
$A(s)+P(s)\se B(s)$ if and only if
  {\em (\ref{eqintfihr1})} holds and
\begin{equation}\label{ceqGg}
 G\leq \sum_{i=1}^{\rho}\min\{r_i, s_i\}+ \max\{e, e'\}.
\end{equation}

\item  
\label{ccolequalcdrns}
If $\br=\bs$, $\br'\neq \bs'$,
let
$$\bar x=\min\{i \; :\; r'_i\neq s'_i\},$$
$$
\bar e=\min\{i\geq\bar x-1\; : \; s'_{i+1}\geq r'_{i+1}\},\quad 
\bar e'=\min\{i\geq\bar x-1\; : \; r'_{i+1}\geq s'_{i+1}\},
$$
$$
\bar G=\rho-1-\sum_{i=1}^{\rho-1}\deg(\gcd(\phi_{i+1}(s, t),\psi_{i+1}(s, t)))- \sum_{i=1}^{\rho}r_i.
$$
Then, there exists a pencil $P(s)\in \FF[s]^{p \times q}$ of $\rank (P(s))=1$ such that 
$A(s)+P(s)\se B(s)$ if and only if
  {\em (\ref{eqintfihr1})} holds and
\begin{equation*}\label{ceqbarGg}
 \bar G\leq \sum_{i=1}^{\rho}\min\{r'_i, s'_i\}+ \max\{\bar e,\bar e'\}.
\end{equation*}

\item \label{cinequal}
If $\br\neq \bs$, $\br'\neq \bs'$,
then there exists a pencil $P(s)\in \FF[s]^{p \times q}$ of $\rank (P(s))=1$ such that 
$A(s)+P(s)\se B(s)$ if and only if   {\em (\ref{eqintfihr1})}  
and one of the four following conditions hold:
\begin{enumerate}
\item[(a)] 
\begin{equation}\label{ccdsr}
\bs\cpr\br, \quad \bs'\cpr\br',
\end{equation}
and
\begin{equation}\label{cbocacdsr}
  \sum_{i=1}^\rho\deg(\lcm(\phi_i(s,t), \psi_i(s,t)))\leq  
 x
  \leq \sum_{i=1}^\rho\deg(\gcd(\phi_{i+1}(s,t), \psi_{i+1}(s,t))),
\end{equation}
where $x= \rho-\sum_{i=1}^{\rho'}r_i-\sum_{i=1}^{\rho'}s'_i$.

\item[(b)]
\begin{equation}\label{cdcrs}
\br\cpr\bs, \quad \br'\cpr\bs',
\end{equation} and
\begin{equation}\label{cbocacdrs}
  \sum_{i=1}^\rho\deg(\lcm(\phi_i(s,t), \psi_i(s,t)))\leq  
 y
  \leq \sum_{i=1}^\rho\deg(\gcd(\phi_{i+1}(s,t), \psi_{i+1}(s,t))),
\end{equation}
where $y=\rho-\sum_{i=1}^{\rho'}s_i-\sum_{i=1}^{\rho'}r'_i$.
\item[(c)]
   {\em (\ref{ccdsr})} and  {\em (\ref{cbocacdrs})} hold.
 \item[(d)]
   {\em (\ref{cdcrs})} and  {\em (\ref{cbocacdsr})} hold.

\end{enumerate}

\end{enumerate}

\end{theorem}

{\it Proof}. It is a consequence of 
Theorem \ref{maintheogenpenc}, Lemmas  \ref{lemmagh2}
 and \ref{lemmapartmin} and Proposition \ref{propconj}.
\hfill $\Box$

\section{Bounds}
\label{seccotas}

Given a pencil $A(s)$ and a  perturbation of it,  $A(s)+P(s)$, where $P(s)$ is a pencil of rank one, and given $\lambda \in \overline{\FF}$,  in this section we obtain bounds for the differences between the generalized Weyr characteristics of $A(s)$  and $A(s)+P(s)$ at $\lambda$.
We include some technical lemmas in Subsection \ref{subsecboundstechnical} and prove the main result in  Subsection \ref{subsecboundsmain}.

The notation used  in this section corresponds to that introduced in Theorem \ref{cmaintheogenpenc}. In particular, when $\br\neq \bs$ and $\br'= \bs'$, the values of $x, e$ and $e'$ are defined in (\ref{T4.8-x}) and (\ref{T4.8-ee'}) (they also appear in Lemma 4.6). Notice that either $e=x-1<e'$ or $e'=x-1<e$, therefore $e\neq e'$.

\subsection{Technical Lemmas} \label{subsecboundstechnical}

\begin{lemma}
\label{lemmareqneqs}
Assume that $\br\neq \bs$, $\br'=\bs'$ and  {\em (\ref{ceqGg})} holds.

\begin{enumerate}
  \item
Case  $e>e'$. 
  \begin{enumerate}
    \item Let $i \in\{x, \dots, e\}$, then
\begin{equation}\label{upe}
 -x-1\leq s_i-r_i\leq -1. \end{equation}
Moreover, if $s_i-r_i=-x-1$, then
\begin{equation}
\label{eqdivphi}
\phi_j(s,t)\mid \psi_j(s,t), \quad 1\leq j \leq \rho.
\end{equation}

\item Let $i > e$, then 
\begin{equation}\label{frome}
 -x\leq s_i-r_i\leq e+1.
\end{equation}
Moreover, if $s_i-r_i=e+1$,  then  
\begin{equation}
\label{eqdivpsi}
\psi_j(s,t)\mid \phi_j(s,t), \quad 1\leq j \leq \rho.
\end{equation}
  \end{enumerate}
\item Case $e'>e$.
  \begin{enumerate}
    \item Let $i \in\{x, \dots, e'\}$, then
\begin{equation*}
 -x-1\leq r_i-s_i\leq -1. \end{equation*}
If $r_i-s_i=-x-1$, then  {\em (\ref{eqdivpsi})} holds.
\item Let $i>e'$, then 
\begin{equation*}
 -x\leq r_i-s_i\leq e'+1.
\end{equation*}
If $r_i-s_i=e'+1$,  then   {\em (\ref{eqdivphi})} holds.
  \end{enumerate}
 
\end{enumerate}

\end{lemma}

{\it Proof.} 
If $\br\neq \bs$ and  $\br'=\bs'$, then
$\rank A(s)=\rank B(s)=\rho$, 
$r_0=s_0=q-\rho$, and
$\min\{e,e'\}=x-1<\max\{e,e'\}$.

As 
$\br \neq \bs$, we have  $\sum_{i=1}^{\rho}r_i\neq 0$ or $\sum_{i=1}^{\rho}s_i\neq 0$, hence $\phi_{1}(s, t)=1$ or
$\psi_{1}(s, t)=1$,  $\deg(\gcd(\phi_{1}(s, t),\psi_{1}(s, t)))=1$, and
$$\sum_{i=1}^{\rho-1}\deg(\gcd(\phi_{i+1}(s, t),\psi_{i+1}(s, t)))=\sum_{i=1}^{\rho}\deg(\gcd(\phi_{i}(s, t),\psi_{i}(s, t))).$$
Bearing in mind that  $\br'=\bs'$ and 
$$
\sum_{i=1}^{\rho}r_i+\sum_{i=1}^{\rho}r'_i+ \sum_{i=1}^{\rho}\deg(\phi_{i}(s, t))=
\sum_{i=1}^{\rho}s_i+\sum_{i=1}^{\rho}s'_i+ \sum_{i=1}^{\rho}\deg(\psi_{i}(s, t))=\rho, $$
we have
$$
G=\sum_{i=1}^{\rho}r_i+\sum_{i=1}^{\rho}X_i-1
=\sum_{i=1}^{\rho}s_i+\sum_{i=1}^{\rho}Y_i-1,
$$
where, for $1\leq i \leq \rho$,
$X_i=\deg (\phi_i(s, t))-\deg(\gcd(\phi_{i}(s, t),\psi_{i}(s, t)))$
and $Y_i=\deg (\psi_i(s, t))-\deg(\gcd(\phi_{i}(s, t),\psi_{i}(s, t))).$
Hence condition (\ref{ceqGg}) is equivalent to
\begin{equation}\label{eqcondA}
\sum_{i=1}^{\rho}(r_i-\min\{r_i, s_i\})+\sum_{i=1}^{\rho}X_i\leq \max\{e, e'\}+1,
\end{equation}
and to
\begin{equation}\label{eqcondB}
\sum_{i=1}^{\rho}(s_i-\min\{r_i, s_i\})+\sum_{i=1}^{\rho}Y_i\leq \max\{e, e'\}+1.
\end{equation}

\begin{enumerate}
\item
Assume that $e>e'$. Then, $e'=x-1<e= \max\{e, e'\} \leq \rho$.
\begin{enumerate}
  \item
Let $i \in\{x, \dots, e\}$. Then $r_i>s_i$ and the upper bound of (\ref{upe}) holds. Moreover,  
\begin{equation*} \label{eqa}
\begin{array}{rl}
\sum_{j=1}^i(r_j-s_j)& =\sum_{j=1}^{e}(r_j-s_j)-\sum_{j=i+1}^{e}(r_j-s_j)\\
& \leq \sum_{j=1}^{e}(r_j-s_j)-(e-i) \\
& \leq\sum_{j=1}^{\rho}(r_j-\min\{r_j, s_j\}) -(e-i).
\end{array}
\end{equation*}
From  (\ref{eqcondA}) we derive
\begin{equation}
\label{eqb}
\begin{array}{rl}
\sum_{j=1}^{i}(r_j-s_j)& \leq 
\sum_{j=1}^{\rho}(r_j-\min\{r_j, s_j\})+\sum_{j=1}^{\rho}X_j-(e-i)\\ 
&\leq (e+1)-(e-i)=i+1,
\end{array}
\end{equation}
then
$$
r_i-s_i=\sum_{j=1}^{i}(r_j-s_j)-\sum_{j=x}^{i-1}(r_j-s_j)\leq (i+1)-(i-x)=x+1,
$$
and the lower bound of  (\ref{upe}) holds.

In the case that $s_i-r_i=-x-1$, from (\ref{eqb}),
 $$i+1\geq \sum_{j=1}^i(r_j-s_j)=\sum_{j=x}^{i-1}(r_j-s_j)+x+1\geq (i-x)+x+1=i+1.$$

Then, again from (\ref{eqb}), we have
$$
(i+1)+ \sum_{j=i+1}^{e}(r_j-s_j)+\sum_{j=e+1}^{\rho}(r_j-\min\{r_j, s_j\})+\sum_{j=1}^{\rho}X_j-(e-i) \leq i+1,$$
i.e.,
$$
\sum_{j=i+1}^{e}(r_j-s_j)-(e-i)+\sum_{j=e+1}^{\rho}(r_j-\min\{r_j, s_j\})+\sum_{j=1}^{\rho}X_j=0.
$$
From this equation we conclude that
$X_j=0,  \ 1\leq j\leq \rho$, which is equivalent to condition
(\ref{eqdivphi}).

\item
  Let $i \in\{e+1, \dots, \rho\}$. Then, from (\ref{eqcondB}),
  \begin{equation}\label{eqsre+1}
s_i-r_i\leq s_i-\min\{r_i, s_i\}\leq \sum_{j=1}^{\rho}(s_j-\min\{r_j, s_j\})+\sum_{i=1}^{\rho}Y_j\leq e+1,
\end{equation}
and, from (\ref{eqcondA}),
\begin{equation*}\label{eqrse+1b}
\begin{array}{rl}
r_i-s_i\leq& r_i-\min\{r_i, s_i\}\leq \sum_{j=e+1}^{\rho}(r_j-\min\{r_j, s_j\})\\
=&\sum_{j=1}^{\rho}(r_j-\min\{r_j, s_j\})-\sum_{j=x}^{e}(r_j-s_j)\\\leq &
\sum_{j=1}^{\rho}(r_j-\min\{r_j, s_j\})+\sum_{i=1}^{\rho}X_j-(e-x+1)\\\leq &(e+1)-(e-x+1)=x.
\end{array}
\end{equation*}
Therefore (\ref{frome}) holds.

In the case that $s_i-r_i= e+1$,  from (\ref{eqsre+1}) we  obtain 
that $Y_j=0,  \ 1\leq j\leq \rho$,
which is equivalent to (\ref{eqdivpsi}).

\end{enumerate}

\item
 The proof is analogous to that of Case 1.

\end{enumerate}

\hfill $\Box$

\bigskip

As a consequence  of Lemma \ref{lemmareqneqs} we obtain the following result.

\begin{lemma}
\label{lemmareqneqs3}
Assume that $\br\neq \bs$, $\br'=\bs'$, and  {\em (\ref{ceqGg})} holds.
Let $a'=\min\{e, e'\}=x-1$ and  $b'=\max\{e, e'\}$. Then
$$
\begin{array}{cl}
  s_i-r_i=0, & 0\leq i \leq a',\\
  -(a'+2)\leq s_i-r_i\leq a'+2, & a'+1\leq i\leq b',\\
  -(b'+1)\leq s_i-r_i\leq b'+1, & i \geq b'+1.
\end{array}
$$
Additionally:

If $s_i-r_i=-(a'+2)$ for some $i\in \{a'+1, \dots, b'\}$ or $s_i-r_i=-(b'+1)$ for some $i>b'$, then  {\em (\ref{eqdivphi})} holds.

If $s_i-r_i=(a'+2)$ for some $i\in \{a'+1, \dots, b'\}$ or $s_i-r_i=b'+1$ for some $i>b'$, then  {\em (\ref{eqdivpsi})} holds.

\end{lemma}

\hfill $\Box$
\begin{lemma}
\label{lemmarandifs1} \
\begin{enumerate}
  \item
Assume that $\rank (A(s))=\rho$ and $\rank (B(s))=\rho+1$,  and that 
 {\em (\ref{ccdsr})} holds. Let $g=\max\{i: r_i>s_i\}$ and 
 $i\in \{g+1, \dots, \rho+1\}$. 
\begin{enumerate}
\item
  If  {\em (\ref{cbocacdrs})} holds, then $0\leq s_i-r_i\leq g.$
  In addition, if $s_i-r_i= g$, then  {\em (\ref{eqdivpsi})} is satisfied. 
\item
  If  {\em (\ref{cbocacdsr})} holds, then $0\leq s_i-r_i\leq g+1.$
  In addition, if $s_i-r_i= g+1$, then   {\em (\ref{eqdivpsi})} is satisfied. 
 \end{enumerate}
 \item 
Assume that $\rank (A(s))=\rho+1$ and $\rank (B(s))=\rho$, and that 
 {\em (\ref{cdcrs})} holds. Let $g=\max\{i: s_i>r_i\}$ and 
 $i\in \{g+1, \dots, \rho+1\}$. 
\begin{enumerate}
\item
  If  {\em (\ref{cbocacdsr})} holds, then $0\leq r_i-s_i\leq g.$
  In addition, if $r_i-s_i= g$, then  {\em (\ref{eqdivphi})} is satisfied.
\item
  If  {\em (\ref{cbocacdrs})} holds, then $0\leq r_i-s_i\leq g+1.$
  In addition, if $r_i-s_i= g+1$, then   {\em (\ref{eqdivphi})} is satisfied. 
 
 \end{enumerate}

\end{enumerate}
\end{lemma}
{\it Proof.}

\begin{enumerate}
  \item
Notice that, from the definition of $g$,  we have  $s_i-r_i\geq 0$.
\begin{enumerate}
\item
From  (\ref{cbocacdrs}) we derive that
$$
\sum_{j=1}^\rho\deg(\lcm(\phi_j(s,t), \psi_j(s,t)))\leq \sum_{j=1}^\rho\deg(\phi_j(s,t))+\sum_{j=1}^{\rho+1}(r_j-s_j).
$$
As a consequence, 
$$
s_i-r_i\leq \sum_{j=g+1}^{\rho+1}(s_j-r_j)\leq \sum_{j=1}^{g}(r_j-s_j)+\sum_{j=1}^\rho X'_j\leq \sum_{j=1}^{g}(r_j-s_j)=g,
$$
where for $1\leq j \leq \rho$, $X'_j=
\deg(\phi_j(s,t))-\deg(\lcm(\phi_j(s,t), \psi_j(s,t)))$.

If $s_i-r_i=g$, then 
$\sum_{j=1}^\rho X'_j=0$,  and (\ref{eqdivpsi})  is satified.
\item
From  (\ref{cbocacdsr}) we derive that
\begin{equation}\label{eqpsi}
\sum_{j=1}^{\rho+1}\deg(\psi_j(s,t))+\sum_{j=1}^{\rho+1}(s_j-r_j)-1\leq \sum_{i=1}^\rho\deg(\gcd(\phi_{i+1}(s,t), \psi_{i+1}(s,t))),
\end{equation}
\begin{itemize}
  \item
    If $\deg(\gcd(\phi_{1}(s,t), \psi_{1}(s,t)))=1$, then $s_i=r_i=0$ and     $g=0$.
  \item
    If $\deg(\gcd(\phi_{1}(s,t), \psi_{1}(s,t)))=0$ then, from (\ref{eqpsi}),
    $$
    \begin{array}{rl}
s_i-r_i & \leq \sum_{j=g+1}^{\rho+1}(s_j-r_j)\leq \sum_{j=1}^{g}(r_j-s_j)+\sum_{j=1}^{\rho+1} Y_j+1 \\
& \leq \sum_{j=1}^{g}(r_j-s_j)+1=g+1,
\end{array}
$$
where  $Y_j=
\deg(\gcd(\phi_j(s,t), \psi_j(s,t))-\deg(\psi_j(s,t))$,  $1\leq j \leq \rho+1$.

If $s_i-r_i=g+1$ then 
$\sum_{j=1}^{\rho+1} Y_j=0$,  and (\ref{eqdivpsi})  is satisfied.
  \end{itemize}

\end{enumerate}

\item The proof is analogous to that of Case 1.
\end{enumerate}
\hfill $\Box$

\subsection{Main theorem} \label{subsecboundsmain}

We now prove the main theorem of the paper. Recall that we use the notation of Theorem \ref{cmaintheogenpenc}.

\begin{theorem}
\label{theomain2}
Given a matrix pencil $A(s)\in \FF[s]^{p\times q}$, let $P(s)\in \FF[s]^{p\times q}$ be a  matrix pencil of rank one.
For $\lambda \in \overline{\FF}$, the following statements hold:
 \begin{enumerate}
\item[(i)]
If both pencils $A(s)$ and $A(s)+P(s)$ are regular, then
$$
-1\leq w_i(\lambda, A(s)+P(s))- w_i(\lambda, A(s)) \leq 1, \quad i \geq 1.
$$
\item[(ii)]
If  $A(s)$ is regular and $A(s)+P(s)$ is singular,  then, taking $a=d_1+1$, we have
$$\begin{array}{cl}
-1\leq  w_i(\lambda, A(s)+P(s))- w_i(\lambda, A(s))\leq 1, & 1\leq i\leq a. \\
-2\leq  w_i(\lambda, A(s)+P(s))- w_i(\lambda, A(s))\leq 0, & i \geq a+1.
\end{array}
$$
\item[(iii)]
If  $A(s)$ is singular and $A(s)+P(s)$ is regular,  then, taking $a=c_1+1$,  we have
$$\begin{array}{cl}
-1\leq  w_1(\lambda, A(s)+P(s))- w_1(\lambda, A(s))\leq 1, & 1\leq i\leq a.\\
0\leq  w_i(\lambda, A(s)+P(s))- w_i(\lambda, A(s))\leq 2,   & i \geq a+1.
\end{array}
$$

\item[(iv)] 
  If both $A(s)$ and $A(s)+P(s)$ are singular:
  \begin{itemize}
  \item Assume that $\rank(A(s))=\rank(A(s)+P(s))$.
  \begin{itemize}
    \item
If $A(s)$ and $A(s)+P(s)$ have the same column minimal indices, then
$$
-1\leq  w_i(\lambda, A(s)+P(s))- w_i(\lambda, A(s))\leq 1,\quad i \geq 1.
$$
\item
  If $A(s)$ and $A(s)+P(s)$ have different column minimal indices, then, taking $a=x$ and $b=\max\{e,e'\}+1$ (notice that $b>a\geq 1$), we have
  \begin{equation} \label{iv1}\hspace{-0.5cm}
\begin{array}{cl}
-1\leq w_i(\lambda, A(s)+P(s))- w_i(\lambda, A(s))\leq 1,
 &\hspace{-0.2cm} 1\leq i \leq a,\\
-(a+1)\leq w_i(\lambda, A(s)+P(s))- w_i(\lambda, A(s))\leq a+1,
   &\hspace{-0.2cm}  a+1\leq i\leq b,\\
-b\leq w_i(\lambda, A(s)+P(s))- w_i(\lambda, A(s))\leq b,
& \hspace{-0.2cm} i \geq b+1.
\end{array}
\end{equation}

  \end{itemize}

\item Assume that $\rank(A(s))<\rank(A(s)+P(s))$. Let  $a=\max\{i: r_i>s_i\}+1$. Then
$$
\begin{array}{ll}
  -1\leq  w_i(\lambda, A(s)+P(s))- w_i(\lambda, A(s))\leq 1, & 1\leq i \leq a,\\
  0\leq  w_i(\lambda, A(s)+P(s))- w_i(\lambda, A(s))\leq a+1, &
i \geq a+1.
  \end{array}
$$
\item Assume that 
  $\rank(A(s))>\rank(A(s)+P(s))$. Let  $a=\max\{i: s_i>r_i\}+1$. Then
  $$
\begin{array}{ll}
  -1\leq  w_i(\lambda, A(s)+P(s))- w_i(\lambda, A(s))\leq 1, & 1\leq i \leq a,\\
  -(a+1)\leq  w_i(\lambda, A(s)+P(s))- w_i(\lambda, A(s))\leq 0, 
& i \geq a+1.
  \end{array}
$$

  \end{itemize}

\end{enumerate}
\end{theorem}
{\it Proof.} We take $B(s)=A(s)+P(s)$. Notice that $-1\leq \rho_2-\rho_1\leq 1$. 
Let
$$
(w^R_1(\lambda, A(s)),  \dots)
=\overline{(n_1(\lambda, A(s)),  \dots, n_{\rho_1}(\lambda, A(s))},
$$$$
(w^R_1(\lambda, B(s)),  \dots)
=\overline{(n_1(\lambda, B(s)),  \dots, n_{\rho_2}(\lambda, B(s))}.$$
By Theorem \ref{carwcanonical}, for $i\geq 1$,
\begin{equation}\label{eqintw0}
w_i(\lambda, B(s))-w_i(\lambda, A(s))=w^R_i(\lambda, B(s))-w^R_i(\lambda, A(s))+s_{i-1}-r_{i-1}.
\end{equation}
By Theorem \ref{cmaintheogenpenc} condition (\ref{eqintfihr1}) holds, then
\begin{equation*}\label{eqintn}
n_{i+\rho_2-\rho_1+1}(\lambda, B(s))\leq n_{i}(\lambda, A(s))\leq n_{i+\rho_2-\rho_1-1}(\lambda, B(s)), \quad i\geq 1, 
\end{equation*}
and, by Lemma  \ref{lemmaconj},  
this is equivalent to
\begin{equation}\label{eqintw}
\rho_2-\rho_1-1\leq w_{i}^R(\lambda, B(s))-w_{i}^R(\lambda, A(s))\leq \rho_2-\rho_1+1, \ \ i \geq 1.
\end{equation}
Additionally, if  (\ref{eqdivpsi}) holds, then 
\begin{equation}\label{eqintw54psi}
\rho_2-\rho_1-1\leq w^R_{i}(\lambda, B(s))-w^R_{i}(\lambda, A(s))\leq \rho_2-\rho_1, \ \ i \geq 1,
\end{equation}
and if (\ref{eqdivphi}) holds, then 
\begin{equation}\label{eqintw54phi}
\rho_2-\rho_1\leq w^R_{i}(\lambda, B(s))-w^R_{i}(\lambda, A(s))\leq \rho_2-\rho_1+1, \ \ i \geq 1.
\end{equation}

We analyze now the different cases enumerated in the statement of the theorem.
\begin{enumerate}

\item[(i)]
If $A(s)$ and $A(s)+P(s)$ are regular, then $p=q=\rho_1=\rho_2$. The result follows from (\ref{eqintw0}) and (\ref{eqintw}).
This result  was proven in \cite[Corollary 4.7]{BaRo18}.
\item[(ii)]
If $A(s)$ is regular and  $A(s)+P(s)$ is singular, then $p=q=\rho_1=\rho_2+1$, $r_0=0$, 
  $(s_1, \dots)=\overline{(d_1)}=(1,\stackrel{(d_1)}{\dots}, 1, 0,\dots)$, and $s_0=1$.
Taking into account  (\ref{eqintw0}) and (\ref{eqintw}), we get
$$
-1\leq  w_i(\lambda, B(s))- w_i(\lambda, A(s))\leq 1,\quad 1\leq i \leq d_1+1,
$$
$$
-2\leq  w_i(\lambda, B(s))- w_i(\lambda, A(s))\leq 0,  \quad i \geq d_1+2.
$$

\item[(iii)]
If $A(s)$ is singular and $A(s)+P(s)$ is regular the proof is analogous to (ii),  exchanging the roles of $A(s)$ and $B(s)$.

  \item[(iv)]
If $A(s)$  and $A(s)+P(s)$ are singular, there are three possibilities: $\rho_1=\rho_2$, $\rho_2=\rho_1+1$ or  $\rho_1=\rho_2+1$.
\begin{itemize}
\item
Let $\rho_1=\rho_2$. As $s_{0}-r_{0}=0$,  from Theorem \ref{cmaintheogenpenc}, necessarily  $\br=\bs$ or $\br'=\bs'$ (notice that if $\bs\cpr\br$, then $r_0=s_0+1$).
\begin{itemize}
\item
If $\br=\bs$, then from conditions  (\ref{eqintw0}) and (\ref{eqintw}) we get
$$
-1\leq  w_i(\lambda, B(s))- w_i(\lambda, A(s))\leq 1,\quad i\geq 1.
$$

\item
  If $\br\neq \bs$, then $\br'= \bs'$ and (\ref{ceqGg}) holds. By Lemma \ref{lemmareqneqs3}, we have
\begin{equation}\label{5.2conseq}
\begin{array}{cl}
  s_{i-1}-r_{i-1}=0 &  1\leq i \leq a,\\
  -(a+1)\leq s_{i-1}-r_{i-1}\leq a+1, & a+1\leq i\leq b,\\
  -b\leq s_{i-1}-r_{i-1}\leq b, &  i\geq b+1.
\end{array}
\end{equation}
Moreover, if (\ref{eqdivphi}) does not hold, then
\begin{equation}\label{5.2no}
\begin{array}{cl}
  -a\leq s_{i-1}-r_{i-1}\leq a+1, & a+1\leq i\leq b,\\
  -b+1\leq s_{i-1}-r_{i-1}\leq b, &  i\geq b+1,
\end{array}
\end{equation}
and, if (\ref{eqdivpsi}) does not hold, then
\begin{equation}\label{5.4no}
\begin{array}{cl}
  -(a+1)\leq s_{i-1}-r_{i-1}\leq a, & a+1\leq i\leq b,\\
  -b\leq s_{i-1}-r_{i-1}\leq b-1, &  i\geq b+1.
  \end{array}
\end{equation}

  We obtain different results depending on the relation between the homogeneous invariant factors of $A(s)$ and $B(s)$.
 
\begin{itemize}
\item If $\phi_i(s,t)=\psi_i(s,t), \ 1\leq i \leq \rho$, then by (\ref{eqintw0}), (\ref{eqintw54psi}), (\ref{eqintw54phi}), and (\ref{5.2conseq}) we obtain
$$
\begin{array}{cl}
  w_{i}(\lambda, B(s))- w_{i}(\lambda, A(s))=0 & \hspace{-0.2cm}  1\leq i \leq a,\\
  -(a+1)\leq w_{i}(\lambda, B(s))- w_{i}(\lambda, A(s))\leq a+1, & \hspace{-0.2cm}  a+1\leq i\leq b,\\
  -b\leq w_{i}(\lambda, B(s))- w_{i}(\lambda, A(s))\leq b, & \hspace{-0.2cm}  i \geq b+1.
\end{array}
$$
\item 
If (\ref{eqdivpsi}) holds and (\ref{eqdivphi}) does not hold, then by (\ref{eqintw0}), (\ref{eqintw54psi}), and (\ref{5.2no}) we obtain
$$
\begin{array}{cl}
-1\leq w_i(\lambda, B(s))- w_i(\lambda, A(s))\leq 0,
 & \hspace{-0.2cm}1\leq i \leq a,\\
-(a+1)\leq w_i(\lambda, B(s))- w_i(\lambda, A(s))\leq a+1,
   & \hspace{-0.2cm}  a+1\leq i\leq b,\\
-b\leq  w_i(\lambda, B(s))- w_i(\lambda, A(s))\leq b,
& \hspace{-0.2cm}  i\geq b+1.
\end{array}
$$
\item
If (\ref{eqdivphi}) holds and (\ref{eqdivpsi}) does not hold, then analogously by (\ref{eqintw0}), (\ref{eqintw54phi}), and (\ref{5.4no}) we obtain
$$
\begin{array}{cl}
0\leq w_i(\lambda, B(s))- w_i(\lambda, A(s))\leq 1,
 &  \hspace{-0.2cm} 1\leq i \leq a,\\
-(a+1)\leq w_i(\lambda, B(s))- w_i(\lambda, A(s))\leq a+1,
   & \hspace{-0.2cm}  a+1\leq i\leq b,\\
-b\leq  w_i(\lambda, B(s))- w_i(\lambda, A(s))\leq b,
& \hspace{-0.2cm}   i \geq b+1.
\end{array}
$$

\item
  If neither (\ref{eqdivphi}) nor  (\ref{eqdivpsi}) are satisfied, then, by (\ref{eqintw0}), (\ref{eqintw}), (\ref{5.2no}), and (\ref{5.4no}), we obtain that  (\ref{iv1}) holds.
	\end{itemize}
\end{itemize}

\item 
  If $\rho_2= \rho_1+1$, 
we have $r_0=s_0+1$ and by Theorem \ref{cmaintheogenpenc}, (\ref{ccdsr}),  and  (\ref{cbocacdsr}) or (\ref{cbocacdrs}), hold.
Let  $g=\max\{i: r_i>s_i\}$, then
$$s_{i-1}-r_{i-1}=-1, \quad  1\leq i \leq g+1.$$

\begin{itemize}
\item
  If (\ref{cbocacdrs}) holds, or (\ref{eqdivpsi}) does not hold, then by Lemma \ref{lemmarandifs1},
  $$
    0\leq s_{i-1}-r_{i-1}\leq g, \quad i\geq g+2.\\
  $$
Therefore, from (\ref{eqintw0}) and (\ref{eqintw}) we obtain
$$
\begin{array}{ll}
  -1\leq  w_i(\lambda, B(s))- w_i(\lambda, A(s))\leq 1, & 1\leq i \leq g+1,\\
  0\leq  w_i(\lambda, B(s))- w_i(\lambda, A(s))\leq g+ 2, 
& i\geq g+2.
  \end{array}
$$
\item 
Alternatively,  if (\ref{cbocacdsr}) and (\ref{eqdivpsi}) hold, then, by Lemma \ref{lemmarandifs1},
$$0\leq s_{i-1}-r_{i-1}\leq g+1, \quad i\geq g+2.$$
Therefore, from (\ref{eqintw0}) and (\ref{eqintw54psi})  we obtain
$$
\begin{array}{ll}
  -1\leq  w_i(\lambda, B(s))- w_i(\lambda, A(s))\leq 0, & 1\leq i \leq g+1,\\
  0\leq  w_i(\lambda, B(s))- w_i(\lambda, A(s))\leq g+ 2, 
& i\geq g+2.
  \end{array}
  $$

\end{itemize}
As $a=g+1$, the result follows.
\item
  If $\rho_1=\rho_2+1$, taking  $g=\max\{i: s_i>r_i\}$, the proof is analogous:

\begin{itemize}
\item 
  If (\ref{cbocacdsr}) holds, or (\ref{eqdivphi}) does not hold, then 
$$
\begin{array}{ll}
  -1\leq  w_i(\lambda, B(s))- w_i(\lambda, A(s))\leq 1, & 1\leq i \leq g+1,\\
  -g-2\leq  w_i(\lambda, B(s))- w_i(\lambda, A(s))\leq 0, 
& i\geq g+2.
  \end{array}
$$
\item 
Alternatively,  if (\ref{cbocacdrs}) and  (\ref{eqdivphi}) hold, then
$$
\begin{array}{ll}
  0\leq  w_i(\lambda, B(s))- w_i(\lambda, A(s))\leq 1, & 1\leq i \leq g+1,\\
  -g-2\leq  w_i(\lambda, B(s))- w_i(\lambda, A(s))\leq 0, 
& i\geq g+2.
  \end{array}
  $$
\end{itemize}
\end{itemize}
\end{enumerate}
\hfill $\Box$

\begin{remark}
We would like to point out that the bounds obtained in Theorem \ref{theomain2} are sharp, i.e., there are examples showing that the bounds are  attained. Concerning the sufficiency,  it is proven in \cite[Corollary 4.14]{BaRo18} that if $A(s)$ and $A(s)+P(s)$ are regular,  the conditions  are sufficient, in the sense that if the bounds are satisfied for some numbers $w_i'$, then there exists a rank-one perturbation $P(s)$ of $A(s)$ such that $w_i(\lambda, A(s)+P(s))=w_i'$. This property is immediately extended to the case where $\br=\bs$.  On the other hand, there are examples showing that in the general case the conditions are not sufficient.
\end{remark}

\section{Conclusion}\label{secconclusion}

We have generalized the notion of Weyr characteristic of an eigenvalue of a pencil (see \cite{Ho90}), and have extended the definition of Jordan chain of square pencils (\cite{LeMaPhTrWiCAOT21}) to arbitrary pencils. Out of them, we have obtained bounds for the changes of the generalized Weyr characteristic of a matrix pencil perturbed by another pencil of rank one. The results in this paper improve the bounds obtained in  \cite[Theorem 7.8]{LeMaPhTrWiCAOT21}. The bounds in Theorem \ref{theomain2} cases $(ii)$ and $(iii)$ are clearly sharper than the corresponding ones in  \cite[Theorem 7.8]{LeMaPhTrWiCAOT21}. Concerning the case $(iv)$, at the cost of splitting the range of indices into different parts, we significantly obtain  better bounds.

It must be remarked that our results hold for any algebraically closed field and for arbitrary rank-one perturbations.

Additionally, we have translated the characterization obtained in \cite[Theorem 5.1]{BaRo20_2} of the changes of the Kronecker structure of a  pencil perturbed by another pencil of rank one, into terms of the  conjugate partitions of the corresponding chains of column and row minimal indices of the pencils involved.

\bibliographystyle{acm}
\bibliography{referencesreg}

\end{document}